\documentclass[11.5pt]{article}

\usepackage{arxiv}
\usepackage{algorithm}
\usepackage{algpseudocode}
\usepackage{caption} 
\usepackage{subcaption}
\usepackage{float}
\usepackage{dblfloatfix}
\usepackage{amsmath}

\usepackage{enumitem}

\newcommand{\R}{\mathbb{R}}

\newcommand{\N}{\mathbb{N}}

\newcommand{\E}{\mathbb{E}}

\newcommand{\Oc}{\mathcal{O}}

\newcommand{\lo}{\longrightarrow}

\newtheorem{remark}{Remark}

\usepackage[utf8]{inputenc} 
\usepackage[T1]{fontenc}    
\usepackage{hyperref}       
\usepackage{url}            
\usepackage{booktabs}       
\usepackage{amsfonts}       
\usepackage{nicefrac}       
\usepackage{microtype}      
\usepackage{lipsum}		
\usepackage{graphicx}
\usepackage{natbib}
\usepackage{doi}
\usepackage{authblk}
\usepackage{orcidlink}

\title{PMBO: Enhancing Black-Box Optimization through Multivariate Polynomial Surrogates}

\author{
Janina Schreiber$^{1, *}$ \orcidlink{0000-0002-8692-0822},
Pau Batlle$^{2, *}$\orcidlink{0000-0003-4886-058X},
Damar Wicaksono$^{1}$\orcidlink{0000-0001-8587-7730} and
Michael Hecht$^{1,3}$\orcidlink{0000-0001-9214-8253}
}

\affil{
$^{1}$ CASUS – Center for Advanced Systems Understanding     Helmholtz-Zentrum Dresden-Rossendorf e.V. (HZDR)\\
 	Untermarkt 20, D-02826 Görlitz, Germany  \\
$^{2}$ California Institute of Technology,	1200 East California Boulevard, Pasadena, CA 91125, US \\
$^3$ University of  Wrocław , Mathematical Institute, pl. Grunwaldzki 2/4 
50-384, Wrocław, Poland\\
$^{*}$ These authors contributed equally to this work}

\date{}

\renewcommand{\headeright}{ }
\renewcommand{\undertitle}{}
\renewcommand{\shorttitle}{Polynomial-Model-Based Optimization}

\hypersetup{
pdftitle={PMBO: Enhancing Black-Box Optimization through Multivariate Polynomial Surrogates},
pdfsubject={q-bio.NC, q-bio.QM},
pdfauthor={Janina Schreiber, Pau Batlle, Damar Wicaksono, Michael Hecht},
pdfkeywords={Polynomial Regression, Bayesian Optimization, Black-box Optimization}
}

\begin{document}
\maketitle

\begin{abstract}
We introduce a surrogate-based black-box optimization method, termed Polynomial-model-based optimization (PMBO). The algorithm alternates polynomial approximation with Bayesian optimization steps, using Gaussian processes to model the error between the objective and its polynomial fit. We describe the algorithmic design of PMBO and compare the results of the performance of PMBO with several optimization methods for a set of analytic test functions.

The results show that PMBO outperforms the classic Bayesian optimization and is robust with respect to the choice of its correlation function family and its hyper-parameter setting, which, on the contrary, need to be carefully tuned in classic Bayesian optimization. Remarkably, PMBO performs comparably with state-of-the-art evolutionary algorithms such as the Covariance Matrix Adaptation -- Evolution Strategy (CMA-ES). This finding suggests that PMBO emerges as the pivotal choice among surrogate-based optimization methods when addressing low-dimensional optimization problems. Hereby, the simple nature of polynomials opens the opportunity for interpretation and analysis of the inferred surrogate model, providing a macroscopic perspective on the landscape of the objective function.
\end{abstract}


\section{Introduction}

Black-box optimization occurs in a wide variety of applications and across all disciplines. 
\cite{Alarie2021} reviews different domains of applications of black-box optimization. These include hyper-parameter selection in \emph{machine learning (ML) algorithms} \cite{DeepLearning, snoek2012practical}, \emph{computer simulations} and even physical or \emph{laboratory experiments} \cite{Alarie2021}.

Scientists and engineers from diverse domains apply black-box optimization methods to problems in computational biology \cite{Ponce_2022} to engineering, such as decision making in autonomous systems \cite{zheng2022gradientfree} and the calibration of hydrological models \cite{Huot_2019}. 

To formalize black-box optimization, we consider a multivariate continuous function $f: \Omega  \to \mathbb{R}$, defined on a bounded domain $\Omega \subseteq [-1, 1]^m$, $m\in \N$  and seek to find a global minimum
\begin{equation}
   x^* =  \mathrm{argmin}_{x \in \Omega} f(x)\,, \quad y^* = f(x^*)\,.
\end{equation}
We assume that no closed-form expression of $f$ is known, and its evaluation is costly but possible in any argument $x \in \Omega$. 

There are two main approaches to black-box optimization:  optimization algorithms \cite{Rios2012}, that evolve a collection of samples in $\Omega$ based solely on the observed function values $f$. 
Popular examples of this class are direct search methods and evolutionary algorithms. Evolutionary algorithms, in particular, like Genetic Algorithms (GA) \cite{GeneticAlgorithm} or Covariance Matrix Adaptation -- Evolutionary Strategy (CMA-ES) \cite{hansen2023cma}, have the advantage that they do not rely on the smoothness or the differentiability of the objective $f$. Whereas GA is inspired by Darwinian evolution and uses recombination and mutation operators, among others, to evolve a population of individuals towards the desired goal, CMA-ES samples from a Gaussian distribution and iteratively adapts the mean and covariance matrix. 

The other, more recent, approach is the use of surrogate models $\hat{f}$ of the function $f$ that guide optimization schemes \cite{Vu2016}. Here, a surrogate approximates the objective function based on the previously observed samples and optimization is guided by the prediction of this surrogate. Bayesian optimization (BO) \cite{BayesianOptimization}, for example, models the objective function as a Gaussian process, iteratively updating its conditional mean and covariance based on the observed samples. Therewith, BO simultaneously delivers a notion of uncertainty for the unseen areas of the landscape, which can be used in an acquisition function to balance the exploitation and exploration rate in the next sample point decision.

Mixture algorithms like the Bayesian Adaptive Direct Search (BADS) \cite{bads} combine the two approaches by switching between evolutionary and model-building methods.

The quality of surrogate-model-based optimization methods goes hand in hand with their approximation power, which has to be strong enough to deliver an accurate prediction of the region of interest, i.e., the region around the optimum. This is of utmost importance when the cost of evaluating $f$ is high. While there is a widespread belief in the infeasibility of using polynomial interpolants $Q_f$ to provide the required approximation power, $\hat f = Q_f \approx f$ we quote Prof Lloyd N. Trefethen who remarked on this persistently held misconception:

\begin{center}
  \fbox{
  \parbox{0.45\textwidth}{\emph{..."one can show that interpolation by polynomials in Chebyshev points is equivalent to interpolation of periodic functions by series of sines and cosines in equispaced points. The latter is the subject of discrete Fourier analysis, and one cannot help noting that whereas there is widespread suspicion that it is not safe to compute with polynomials, nobody worries about the Fast Fourier Transform (FFT)! In the end, that this may be the biggest difference between Fourier and polynomial interpolants, the difference in their reputations."} \citep{trefethen2019}.}}
\end{center}
We underline Trefethen's statement with theoretical results and successful practical applications in the sections below.

\subsection{Contribution}
In this work, we propose a surrogate-model-based optimization algorithm that uses polynomial interpolation to approximate $f$. The algorithm, \emph{Polynomial Model Based Optimization (PMBO)} builds on recent contributions \cite{Hecht2017Quadratic,Hecht2018Interpolation,IEEE,MIP} that extend the aforementioned one-dimensional (1D) polynomial interpolation techniques to multi-dimensions (mD), by resisting the \emph{curse of dimensionality}, reaching exponential approximation rates for a class of functions called \emph{Bos-Levenberg-Trefethen (BLT)} \citep{trefethen2017a, bos2018bernstein}. 

In addition to the polynomial surrogates $\hat f =Q_f$, PMBO follows a \emph{Bayesian optimization} strategy, based on a Gaussian process $Y$ -- informing on its uncertainty -- with $\E[Y(x)] = Q_f(x)$. This is equivalent to using a centered Gaussian process to fit the error between the true function $f$ and the model $Q_f$.
This enables iterative refinement of $Q_f$, balancing exploration and exploitation of the sampling process.
 
As we argue and demonstrate, PMBO delivers an optimization algorithm that tremendously strengthens the surrogate-building step.
Contrary to the classic Bayesian optimization, PMBO appears to be largely independent of the particular choice of the Gaussian process;
PMBO shows high robustness to varying ruggedness of the objective function, tightly and efficiently localizes the desired optima, and 
``on the fly'' approximates the objective function landscape.

\newpage
\section{Approximation theory}
\label{sec:approximation-theory}

We discuss the main aspects of surrogate modelling in the cases of \emph{polynomials and Gaussian processes}.

\subsection{Polynomial approximation theory}

Polynomial interpolation goes back to Newton, Lagrange, and others  \citep{LIP}, and its fundamental importance in mathematics and computing is undisputed.
Interpolation is based on the fact that, in 1D, one and only one polynomial $Q_{f,n}$ of degree $n$ can interpolate a function
$f : \R \lo \R$ in $n+1$ distinct \emph{unisolvent interpolation nodes} $P_n \subseteq \R$, $Q_{f,n}(p_i) =f(p_i)$ for all $p_i \in P_n$, $0\leq i \leq n$.

Though the famous \emph{Stone-Weierstrass approximation theorem} \citep{weier1,weier2}
states that even in multi-dimensions (mD) any continuous function $f : \Omega \to \R$ can be uniformly approximated by
polynomials, this does not necessarily apply for interpolation. In contrast to interpolation, the Weierstrass approximation theorem does not require the polynomials to coincide with $f$ anywhere, meaning there might be a sequence of polynomials $Q_{f,n}$ with $Q_{f,n}(x) \not = f(x)$ for all $x\in \Omega$, but still
\begin{equation}
Q_{f,n} \xrightarrow[n \rightarrow \infty]{} f \quad \text{uniformly on}\,\,\, \Omega\,.
\end{equation}
There are several constructive proofs of the Weierstrass approximation theorem, including the prominent version given by Serge~\citet{bernstein1912}. Although the resulting Bernstein approximation scheme
is \emph{universal} -- approximating any continuous function -- and has been proven to reach the optimal (inverse-linear) approximation rate for the absolute value function $f(x) = |x|$ \citep{bernstein1914},
it achieves only slow convergence rates for analytic functions, resulting in a high computational cost in practice.

In contrast, interpolation in Chebyshev and Legendre nodes is known to be \emph{non-universal} \citep{faber},
but ensures the approximation of the broad class of Lipschitz continuous functions -- Runge's overfitting phenomenon completely disappears --  with fast exponential
approximation rates \citep{trefethen2019}, appearing for analytic functions.

Consequently, there has been much research into mD extensions of 1D interpolation schemes and their approximation capabilities.
The question of how to \emph{resist the curse of dimensionality} is particularly central to these endeavours. 

\subsection{Multivariate interpolation theory}\label{sec:MIT}

Recent results \cite{bos2018bernstein,trefethen2017a,MIP} provide a solution for the class of \emph{Bos-Levenberg-Trefethen (BLT) functions} that can be approximated by interpolation with exponential rates, resting on sub-exponentially many interpolation nodes.
We shortly outline the  essential aspects of this solution.  

In contrary to 1D, there is no canonical notion of \emph{polynomial degree in mD}: We consider the multi-index sets 
\begin{equation}
\label{eq:A}
    A_{m,n,p}=\{\alpha \in \mathbb{N} \mid \| \alpha \|_p \le n\} \subseteq{\mathbb{N}^m\,, m,n \in \N\,, p > 0 }
\end{equation} 
where $\|\cdot\|_p$ denotes the $l_p$-norm and we assume that $A_{m,n,p}$ is lexicographically ordered.
The $A_{m,n,p}$ induce the polynomial spaces $\Pi_{m,n,p}=\mathrm{span}\{ x^\alpha = x_1^{\alpha_1}\cdots x_m^{\alpha_m} \}$ of bounded $l_p$-degree. 

The class of BLT functions is a subclass of continuous functions $C^0(\Omega,\R)$ that can be analytically (holomorphically) extended to the Trefethen domain $N_{m,\rho}$, $\rho >1$ -- a generalized  mD version of a Bernstein ellipse.
Given a BLT function $f$, the following approximation rates 
\begin{equation}\label{eq:TR}
  \| f - Q_{A_{m,n,p}}(f)\|_{C^0(\square_m)}  = \begin{cases}
      \Oc(\rho^{-n/\sqrt{m}})  &\,, \quad p =1 \\
\Oc(\rho^{-n}) &\,, \quad p =2\\
\Oc(\rho^{-n})  &\,, \quad p =\infty
  \end{cases}
\end{equation}
can theoretically be reached at best. 

\begin{remark}
Note that the number of coefficients 
\begin{equation*}
  |A_{m,n,p}| \approx  \begin{cases}
      \binom{m+n}{n} \in \Oc(m^n) &\,, \quad p =1 \\
\frac{(n+1)^m }{\sqrt{\pi m}} (\frac{\pi \mathrm{e}}{2m})^{m/2} \in o(n^m) &\,, \quad p =2\\
(n+1)^m \in \Oc(n^m)  &\,, \quad p =\infty
  \end{cases}
\end{equation*}
for total degree interpolation, $p=1$, scales polynomially,
for Euclidean degree, $p=2$, scales sub-exponentially, whereas
for maximum degree, $p = \infty$, scales exponentially with the dimension $m \in \N$. Thus, due to the exponential rate, Eq.~\eqref{eq:TR}, interpolation of BLT functions with respect to Euclidean degree may resist the curse of dimensionality, while interpolation with total or maximum degree is sub-optimal.

\end{remark}

In fact, the optimal convergence rates Eq.~\eqref{eq:TR} are reached by the multivariate interpolation algorithm (MIP) \cite{MIP} implemented in \cite{minterpy} which is capable for addressing the non-tensorial interpolation task for any $p >0$.  
This enables practical and efficients computations of accurate approximations for regular (BLT) functions.

\subsection{Approximation by Gaussian process regression}
\label{sub:gp}

Gaussian process regression \cite{bishop2006pattern,williams2006gaussian} approximates functions by conditioning a Gaussian process to the observed data.
In particular, given a Gaussian process (GP) with a prior mean function $m$, a correlation kernel function $k$, and a process variance $\sigma^2$ written as
\begin{equation*}
   Y(x) \sim \mathop{GP}(m(x), \sigma^2 k(x, x')),
\end{equation*}
the conditional expectation of the process at a point $x$ after observing the data $\{\mathcal{X}, \mathcal{Y}\}=\{(x_i, y_i)\}_{i = 1}^N$ is given by the function
\begin{equation}
\label{eq:gpr}
  \mu_{Y | \mathcal{Y}} (x) = m(x) + c(x) C^{-1} \left( \mathcal{Y} - m(\mathcal{X}) \right)
\end{equation}
where, for $x_i, x_j \in \mathcal{X}$,
\begin{itemize}
    \item $C = \sigma^2 \left(  k(x_i, x_j) \right)$ is the covariance matrix of $\mathcal{X}$,
    \item $c(x) = \sigma^2 \left( k(x, x_j) \right)$ is the covariance vector between the point $x$ and $\mathcal{X}$, and
    \item $m(\mathcal{X}) = \left( m(x_i) \right)$ is the vector of prior mean values evaluated at $\mathcal{X}$.
\end{itemize}

The Reproducing Kernel Hilbert Space (RKHS) perspective to GP Regression \cite{1807.02582} shows that the centered GP $\mu_{Y | \mathcal{Y}} - m$ in \eqref{eq:gpr} solves the optimization problem
\begin{equation*}
\min_{f} \quad  \|f\|_k \quad 
\text{s.t.}  \quad  f(\mathcal{X}) = \mathcal{Y} - m(\mathcal{X})\,,
\end{equation*}
where $\|\cdot\|_k$ is the norm of the RKHS induced by the kernel $k$.
In other words: $f$ is the most regular function in the RKHS that interpolates the data $\mathcal{X}$ when the GP is centered, i.e., with the mean subtracted.
Hereby, the choice of the kernel function $k$ is implicitly selecting the space of candidate functions, inferring $f$.
Convergence guarantees and rates for Gaussian processes have been studied in different settings,
including classic results \cite{Choi2007}, data-driven covariance functions \cite{Teckentrup2020}, and possibly misspecified settings \cite{wynne2021convergence, wang2022convergence}.

\section{Methodology}
\label{sec:methodology}

PMBO rests on the recent results on polynomial interpolation \cite{MIP} and regression \cite{2212.11706}. Whereas we detail the theoretical aspects in Section~\ref{sec:MIT} and \ref{sub:gp}, we summarize the implementation aspects in section \ref{chap: polynomial interpolation} and \ref{chap: polynomial regression}. 
The methodological ingredients of PMBO as provided in Section~\ref{chap:algorithm}.
\subsection{Multivariate polynomial interpolation (MIP)}
\label{chap: polynomial interpolation}

In 1D, the \emph{barycentric Lagrange interpolation} \citep{berrut,trefethen2019} is known to be the optimal choice.
It enables a numerically stable interpolation up to degree $n = 1.000.000$ or more.
In higher dimensions, however, only polynomial degrees of up to $n = 100$ might be computable at most.
That is why MIP extends 1D Newton interpolation, known to be stable in this range of degrees for Leja ordered nodes \citep{tal1988high,leja}.

Let $A=A_{m,n,p}$ be a lexicographically ordered multi-index set, Eq.~\eqref{eq:A},
$\mathrm{Cheb}_n = \{ \cos(\frac{k\pi}{n}) :  0 \leq k \leq n\}$
be the Chebyshev-Lobatto nodes, which we assume to be Leja ordered,
then the \emph{multivariate Newton polynomials} are given by
\begin{equation}
\label{Newt}
  N_\alpha(x) = \prod_{i=1}^m\prod_{j=0}^{\alpha_i-1}(x_i-q_{j}) \,, \quad  q_j \in \mathrm{Cheb}_n\,, \alpha \in A\,.
\end{equation}

MIP is based on a multi-dimensional extension of the classic Newton-divided difference scheme that computes the Newton interpolant
\begin{equation}\label{eq:QF}
 Q_{f,A}(x) = \sum_{\alpha \in A} c_\alpha N_\alpha(x)
\end{equation}
of any function $f : [-1,1]^m \to \R$ in quadratic runtime $\Oc(|A|^2)$ and linear storage $\Oc(|A|)$.
The evaluation of $Q_{f,A}(x_0)$ in any argument $x_0 \in \R^m$ is stably  realised in $\Oc(m|A|)$,
whereas $k$-th order differentiation requires $\Oc(mn^k|A|)$ \citep{MIP}. 

In particular, MIP enables the computation of the formulas for the \emph{multivariate Lagrange polynomials}
\begin{equation}\label{eq:Lag}
 L_\alpha(x) = \sum_{\beta \in A} c_{\alpha,\beta} N_\beta(x)\,, \quad L_{\alpha}(p_\gamma) = \delta_{\alpha,\gamma}\,, \quad c_\alpha \in \R\,,
\end{equation}
where $p_\gamma = (q_{\gamma_1},\dots,q_{\gamma_m})\,, q_i \in \mathrm{Cheb}_n$
are the Leja ordered Chebyshev-Lobatto (LCL) nodes. This extends Eq.~\eqref{eq:QF} to 
\begin{equation}\label{eq:QF2}
 Q_{f,A}(x) = \sum_{\alpha \in A} f(p_\alpha) L_\alpha(x) = \sum_{\alpha \in A} c_\alpha N_\alpha(x)\,,
\end{equation}
reflecting the Lagrange or Newton interpolation. 
As outlined in Section~\ref{sec:MIT}, the proper choice of $A$ is essential to resist the curse of dimensionality and is set to $A =A_{m,n,2}$ by default. 

The approximation task occurring in the PMBO surrogate modelling, however, requires dealing with unstructured, scattered node sets;
this demands extending MIP to a general fitting scheme. 

\subsection{Multivariate polynomial regression}
\label{chap: polynomial regression}

To provide PMBO with the freedom to sample the objective function at any point apart from the interpolation nodes, we introduce the Newton-Lagrange regression.
That is, we obtain the coefficients of the polynomial model by solving the least squares problem
\begin{equation}\label{eq:REG}
    \hat{c} =\operatorname{argmin}_{c \in \mathbb{R}^{|A|}}\left\|R_A c - F\right\|^2,    
\end{equation}
where $F$ is the vector of observed function values at the sample points $\{x_i\}_{i = 1}^{N}$.
The regression matrix  $R_A$ is obtained by evaluating the Lagrange polynomials (Eq.~\eqref{eq:Lag}) at the sample points
\begin{equation}
\label{eq:RA} 
    R_A=\left((L_\alpha (x_i)\right)_{i=1, \ldots,N, \alpha \in A} \in \mathbb{R}^{N\times|A|}\,.
\end{equation}
We assume $R_A$ to be of full rank $(|A| \leq N)$ such that the optimization problem in Eq.~\eqref{eq:REG} is well-posed.
As the Lagrange basis turns out to be close to an orthogonal basis, our particular choices ensure numerical stability and strong approximation power of the regression scheme \cite{REG_arxiv}. 

Having presented the ingredients of the  polynomial surrogate modelling, we focus next on the algorithmic aspects of PMBO.  


\subsection{Algorithmic design of PMBO}
\label{chap:algorithm}

Similarly to other surrogate-model-based optimization algorithms, our algorithm runs in three steps: 
\begin{enumerate}[label=(\roman*)]
    \item \textbf{Design}: Sample $N_0$ points from the objective function as starting points for the algorithms.
    \item \label{II}\textbf{Model}: Construct a surrogate model with the previous data.
    \item \textbf{Search}: Use the surrogate to decide on the next point to evaluate. If termination conditions are not met, go to step \ref{II} and enhance the surrogate model with a new data point.
\end{enumerate}
In the following we detail the three steps of this methodology.

\subsubsection{Initialization}
\label{subsub:pmbo-initialization}

The algorithm starts by fitting an initial multivariate polynomial surrogate model to the objective function. Therefore, a low polynomial degree ($n=2$ by default) is used, which directly influences the number of coefficients and, hence, the size of the regression matrix. To achieve numerical stability of the least-square regression, Eq.~\eqref{eq:RA}, the number of initial samples $N_0$ is chosen larger than the number of coefficients, $|A| \ll N_0$.

Any sampling method may be used to create the initial samples, e.g., simple random sampling, Sobol' sequence, or LCL nodes, Eq~\eqref{eq:Lag}.

\subsubsection{PMBO surrogate model}
\label{subsub:pmbo-surrogate}
For a given sample, a multivariate polynomial is fitted via least-squares (Section~\ref{chap: polynomial regression}).
The \emph{multivariate interpolation algorithm (MIP)} is implemented in the {\sc minterpy} package \cite{minterpy}, conducting the interpolation and regression tasks appearing as part of the PMBO computations. \\
While MIP operates on the hypercube $[-1, 1]^m$,  re-scaling from the respective domain to $[-1, 1]^m$ is realized via  
the \emph{min-max} transformation 
\begin{equation}
    x_i = -1 + \frac{1 - (-1)}{ub_i - lb_i} (x_{o, i} - lb_i), \; i = 1, \ldots, m\,,
\end{equation} 
where $x_{o, i}$ and $x_i$ are the $i$-th dimensional sample values in $[lb_i, ub_i]$ and $[-1, 1]$; 
$lb_i$ and $ub_i$ are the lower and upper bounds of the $i$-th dimension, respectively.
Optionally, a logarithm may be applied to the above transformation to enable searching the domain with less biased distributions of samples across different orders of magnitude.

Once fitted, the multivariate polynomial surrogate model $Q_f$ serves as the prior mean to the Gaussian process
\begin{equation*}
   Y(x) \sim \mathop{GP}(Q_f(x), \sigma^2 k(x, x'))\,,
\end{equation*}
where $k$ is the chosen correlation kernel function.
Given the sampled (observed) data $\mathcal{X}$, $\mathcal{Y}$, and the polynomial mean $Q_f$,
the posterior (conditional) mean and variance of the Gaussian process at point $x$ are 
\begin{align*}
\mu_{Y | \mathcal{Y}} (x) & = Q_f(x) + c(x) C^{-1} \left( \mathcal{Y} - Q_f(\mathcal{X}) \right) \\
\sigma_{Y | \mathcal{Y}}^2 (x) & = \sigma^2 - c(x) C^{-1} c^T(x) + u^T(x) \left( R_A^T C^{-1} R_A \right) u(x),
\end{align*}
where $c(x)$ and $C$ are as defined in Section~\ref{sub:gp};
$\sigma^2$ is the stochastic process variance (Section~\ref{sub:gp});
$R_A$ is the polynomial surrogate regression matrix (Eq.~(\ref{eq:RA}));
and $u(x) = R_A^T C^{-1} c^T(x) - r_A(x)$ with $r_A(x) = \left( L_{\alpha}(x) \right)_{\alpha \in A}$ (Sections~\ref{chap: polynomial interpolation} and \ref{chap: polynomial regression}).



Both the posterior mean and variance are the main ingredients for the acquisition function that decides on the next sample point.

\subsubsection{Acquisition function}
\label{subsub:pmbo-acquisition} For balancing exploration and exploitation of the search the acquisition functions is maximized, e.g., by the Broyden-Fletcher-Goldfarb-Shanno (BFGS) algorithm. Potential acquisition functions are
\begin{itemize}
    \item Expected Improvement \cite{Jones1998},
    \item Probability of Improvement \cite{Kushner1964ANM}, or 
    \item Upper Confidence Bound \cite{Auer2003UsingCB}.
\end{itemize}
For the Gaussian process as well as the acquisition function, we make use of the implementation provided by the Python package {\sc BayesO} \cite{Kim_BayesO_A_Bayesian_2023}.
\subsubsection{Updating the model complexity}
\label{subsub:pmbo-update}

As mentioned in Section \ref{subsub:pmbo-initialization}, the polynomial is initialized with a low-degree polynomial, which may not be sufficient to capture the complexity of the objective function.
A low-degree polynomial might be able to approximate the macro-landscape of the objective function as, for instance, in the case of the \textit{Rastrigin} (refer to Figure~\ref{fig:Rastrigin-landscape}) function.
However, when considering, e.g., the \textit{Schwefel}  (refer to Figure~\ref{fig:Schwefel-landscape})  function with its several local minima, a polynomial of degree $n=2$ or $n=3$ cannot capture the global structure of the function. For resolving the issue,  
we increase the polynomial degree, whenever enough samples have been evaluated, meaning that the sample size exceeds the number of coefficients of the polynomial.
This condition is posed to fulfill the assumption that $R_A$ is of full rank (Section \ref{chap: polynomial regression}).

The polynomial surrogate is iteratively updated by either adding a new sample to the regression matrix or by increasing the polynomial degree. The updated polynomial serves then as the mean for the Gaussian process. It is worth noting that the Gaussian process is re-defined in every iteration using the new constructed surrogate as the mean, instead of updating the mean in a Bayesian way. This is a key distinction between our method and any type of Bayesian optimization. \\
 A formal combination of the ingredients of Sections~\ref{subsub:pmbo-initialization} to \ref{subsub:pmbo-update}, describing the workflow of PMBO, is given in Alg.~\ref{alg:PMBO}.

\begin{algorithm}[H]
\caption{PMBO}
\begin{algorithmic}
\State $m \gets \mathrm{dimension}$
\State $n \gets \mathrm{polynomial \; degree}$
\State $N \gets \mathrm{max. \; number \; of \; objective \; function \; evaluations}$
\State $N_0 \gets \mathrm{number \; initial \; samples}$
\State $blackbox \gets \mathrm{objective \; function}$
\State $x \gets get\_initial\_samples(N_0)$
\State $y \gets blackbox(x)$
\If{$N_0 \le N$}
    \State $N_0 \gets N_0+1$
    \State $x_t \gets transform(x)$
    \If{$nr\_of\_coefficients(n+1) < length(x)$}
        \State $poly.increase\_poly\_degree(n+1)$
    \EndIf
    \State $poly \gets fit\_polynomial(x_t, y)$
    \State $x \gets transform\_back(x_t)$
    \State $x_{predict} \gets generate\_test\_samples()$
    \State $prior\_mean \gets poly.predict(x_{predict})$
    \State $gaussian\_dist \gets gp.predict(x_t, y, x_{predict}, prior\_mean)$
    \State $x_{next} \gets acquisition\_fct.choose\_next\_sample(gaussian\_dist)$
    \State $y_{next} \gets blackbox(x)$
    \State $x \gets append(x, x_{next})$
    \State $y \gets append(y, y_{next})$
    \State $y_{best} \gets min(y)$
\EndIf    
\end{algorithmic}
\label{alg:PMBO}
\end{algorithm}

\section{Results}
\label{sec:results}
For benchmarking, we use a set of analytic functions of the black-box optimization benchmarking (BBOB) test suite, in the COCO \cite{hansen2021coco} standard test suite that contains 24 noiseless, single-objective and scalable test functions. 
We use the implementation provided by the package {\sc cma} \cite{cmabbobbenchmarks} and consider five test functions -- adjustable for different dimensions -- each of which is indexed by an ID as encoded in Table~\ref{tab:functions} following the convention of BBOB \cite{Hansen2010}.
\begin{center}
\begin{table}[h]
\centering
\caption{Each of the analytic benchmark functions is indexed by an ID, the function name is specified in the according column and the minimum is denoted in the column $f(x^*)$.}
\begin{tabular}{||c c c||} 
 \hline
 ID & name & $f(x^*)$\\ [0.5ex] 
 \hline\hline
  $1$ & Sphere &  $-92.65$ \\ 
 \hline
  $6$ & Attractive Sector & $83.48$\\
 \hline
  $10$ & Ellipsoidal & $-78.99$\\
 \hline
  $15$ & Rastrigin  & $-44.77$ \\
 \hline
 $20$ & Schwefel & $183.12$ \\ 
 \hline
\end{tabular}

\label{tab:functions}
\end{table}
\end{center}
We compare PMBO with the following three different methods along with the selected respective implementation in Python:
\begin{itemize}
    \item Classic Bayesian optimization \cite{BayesianOptimization} as implemented in the package {\sc BayesO} \cite{bads},
    \item Covariance Matrix Adaptation - Evolution Strategy (CMA-ES) \cite{hansen2023cma} as implemented in the package {\sc cma} \cite{cmabbobbenchmarks}, and
    \item Bayesian Adaptive Direct Search (BADS) \cite{bads} as implemented in the package {\sc PyBads} \cite{bads}.
\end{itemize}
Unless explicitly stated otherwise in the subsequent sections, we use the default settings for each of the packages above.
We carry out numerical experiments to optimize the functions in Table~\ref{tab:functions} using PMBO and the three other aforementioned methods. Therefore, a budget of objective function evaluations is fixed to $N=m \cdot 100$, within which the lowest objective function value serves as our performance metric.

Due to the random nature of all the considered methods, we replicate each optimization run $5$ times and statistically summarize the results.

\subsection{PMBO vs. Bayesian optimization}
\label{sec:PMBOvsBO}

One of the key challenges when fitting Gaussian processes is to find an appropriate set of hyper-parameters,
e.g., a suitable kernel, the corresponding range parameter $l$.
In the setting of black-box optimization, it might be hard to a-priori know which set of parameters is suitable for the underlying objective function and hence different strategies for overcoming this issue are common in the field; examples are \textit{Bayesian model selection} \cite{sanzalonso2023inverse, garnett_2023} or \textit{model averaging} \cite{garnett_2023}.

To analyse how PMBO compares to the classic Bayesian optimization, we fix the kernel, as well as the range parameter $l$ and the stochastic process variance and examine the performance of both algorithms for every combination of those \cite{damar2018}.
 We run experiments for dimension $m=2$ and $m=5$ with the kernel methods Matern$_{3/2}$, Matern$_{5/2}$ and Squared Exponential.
 The range and stochastic variance parameter take the values $l = \sigma^2 = \{ 10^{-3}, 10^{-2}, 10^{-1}, 1, 10, 100,\\  500, 1000 \}$ for dimension $m=2$ and are set to  $l = \sigma^2 = \{ 10^{-3},  1, 1000 \}$ for dimension $m=5$.

In the following, we term Bayesian optimization with fixed hyper-parameters as $BO_{\mathrm{fixed}}$.
We use an in-house implementation based on the Python package {\sc minterpy} for fitting the multivariate polynomial \cite{casusmin47:online} and the Python package {\sc BayesO} \cite{Kim_BayesO_A_Bayesian_2023} for running the Bayesian optimization process.
All other parameters --- besides those mentioned above --- are kept at their default value set by {\sc BayesO} \cite{Kim_BayesO_A_Bayesian_2023}.


We choose \emph{simple random sampling} as the initial sampling strategy (Section~\ref{subsub:pmbo-initialization}) and \emph{expected improvement} as the acquisition function (Section~\ref{subsub:pmbo-acquisition}).
To ensure a fair comparison, the initial samples are identical for both methods in each replicate.

\begin{figure*}[!hb]
\centering
\begin{subfigure}[b]{0.9\linewidth}
    \includegraphics[width=1.\textwidth]
    {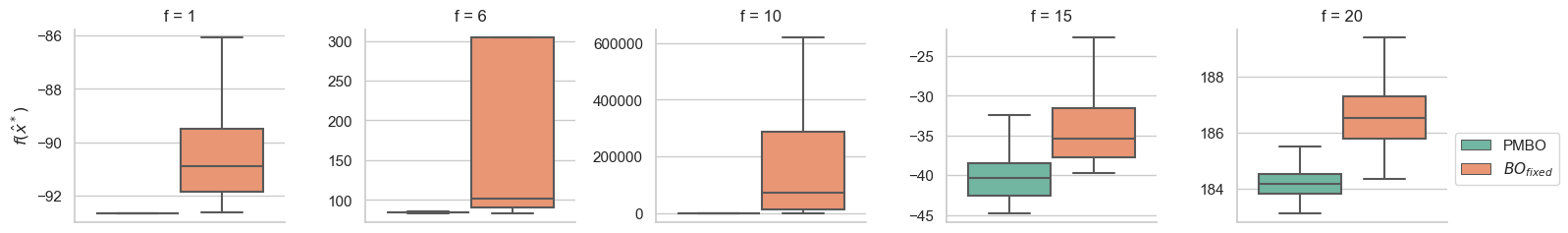}
    \caption{Dimension $m=2$}
    \label{fig:combined_d2}
\end{subfigure}
\begin{subfigure}[b]{0.9\linewidth}
    \includegraphics[width=\textwidth]{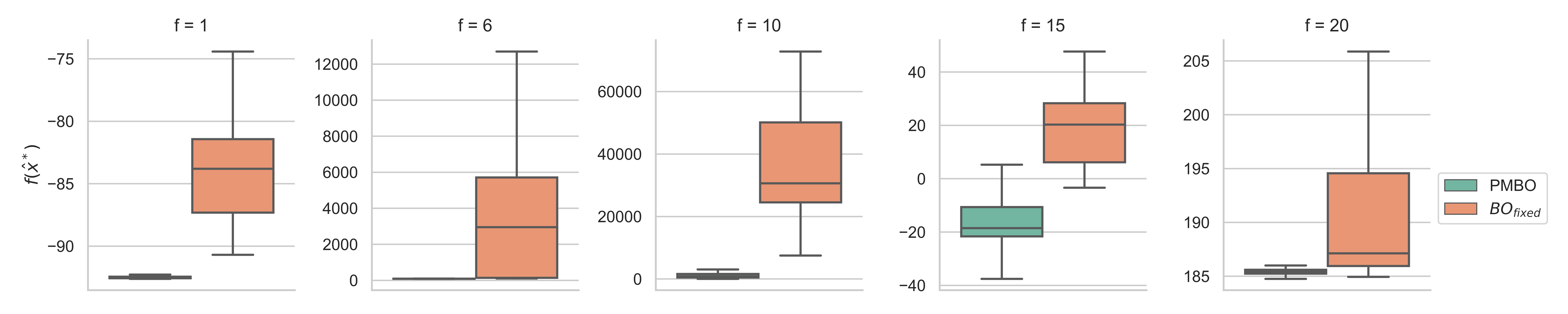}
    \caption{Dimension $m=5$}
    \label{fig:combined_d5}
\end{subfigure}
\caption{Combined results over all possible kernel and hyper-parameter combinations for dimension $m=2$ and dimension $m=5$.}
\label{fig:combined}
\end{figure*}

Figure~\ref{fig:combined} presents the results for the two algorithms over all considered kernel and hyper-parameter combinations for dimension $m=2$ and $m=5$.  Each column represents one of the respective functions.

Clearly, PMBO outperforms $BO_\mathrm{{fixed}}$ in the examined showcases.
Not only the median of the predicted minimum is closer to the true optimum, but also its distribution (over all considered kernels and hyper-parameters as well as replications) is much more narrow, indicating that PMBO is less sensitive to the considered hyper-parameters.

\begin{figure}[h]
\begin{center}
\begin{subfigure}[b]{.49\linewidth}
\includegraphics[width=1.\linewidth]{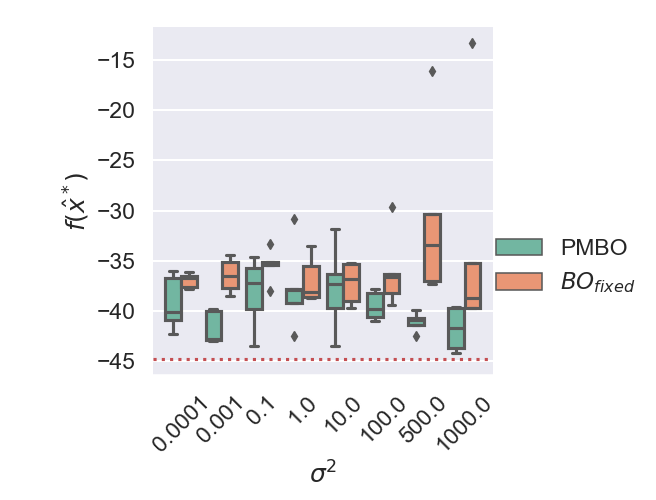} 
\caption{Lowest objective function value $f(\hat{x}^*)$ for the \emph{Rastrigin} function $(ID=15)$ and Bayesian optimization using the Matern$_{5/2}$ kernel and range parameter $l=100$.}
\label{fig:signalf15}
\end{subfigure}
\begin{subfigure}[b]{.49\linewidth}
\includegraphics[width=1.\linewidth]{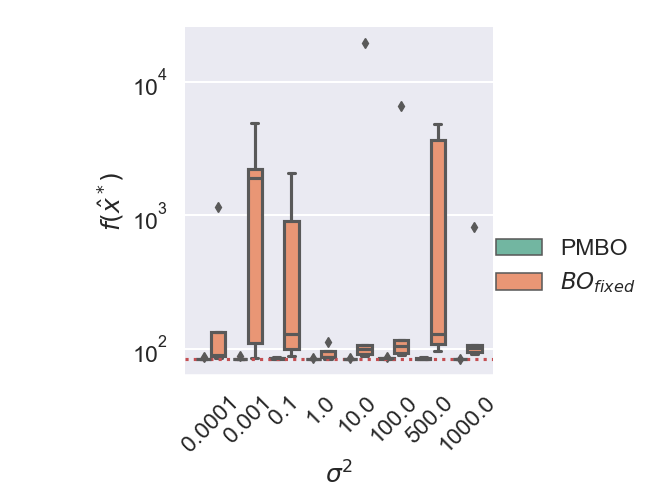} 
\caption{Lowest objective function value $f(\hat{x}^*)$ for the \emph{Attractive Sector} ($ID=6$) and Bayesian Optimization using the Matern$_{3/2}$ kernel and a fixed range parameter $l=0.1$.}
\label{fig:signalf6}
\end{subfigure}
\caption{Lowest objective function values from Bayesian optimization using different fixed choices of kernels and range parameter over varied $\sigma^2$.
The red dotted lines in both figures denote the respective true minimum.}
\label{fig:fixedparams}
\end{center}
\end{figure}

The box plot in Figure~\ref{fig:fixedparams} is a concrete example that underlines the statement above: in the case of Figure~\ref{fig:signalf15}, PMBO is more robust towards changes of $\sigma^2$ than $BO_{\mathrm{fixed}}$ on the \textit{Rastrigin} function ($ID = 15$).
For the function \textit{Attractive Sector} ($ID = 6$) in Figure~\ref{fig:signalf6} the difference exaggerates even more.
Here, a poor choice of $\sigma^2$ drastically alters the performance of $BO_{\mathrm{fixed}}$, whereas PMBO performs relatively stable across a wide range of the hyper-parameter values.

\begin{figure}[H]
\centering
\begin{subfigure}[b]{0.32\linewidth}
    \includegraphics[width=1.\textwidth]
    {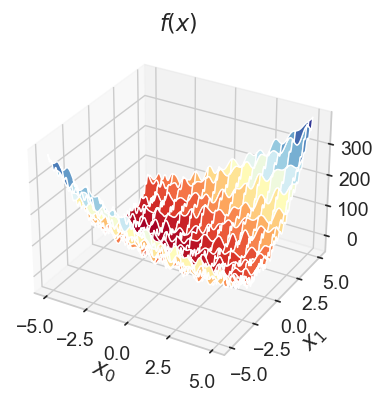}
    \caption{Rastrigin function}
    \label{fig:Rastrigin-landscape}
\end{subfigure}
\begin{subfigure}[b]{0.35\linewidth}
    \includegraphics[width=1.\textwidth]
    {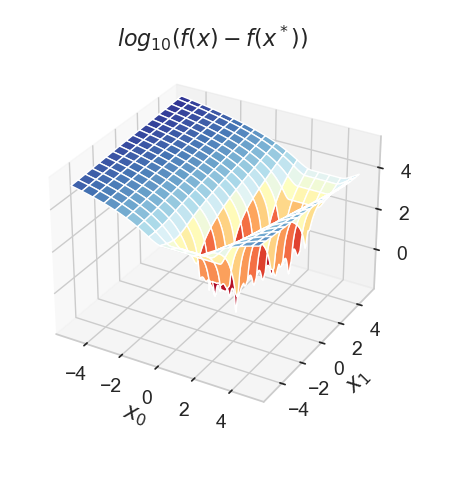}
    \caption{Schwefel function}
    \label{fig:Schwefel-landscape}
\end{subfigure}
\begin{subfigure}[b]{0.8\linewidth}
    \includegraphics[width=1\textwidth]{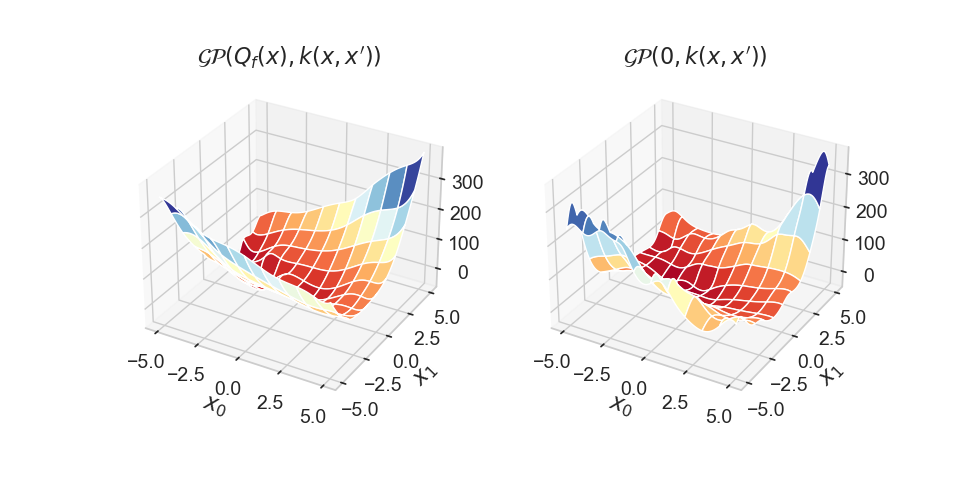},
    \caption{The prediction of the Gaussian mean after seeing the polynomial as prior $\mathop{GP}(Q_f(x), k(x, x')$ (left) compared to the prediction of a Gaussian with zero mean $\mathop{GP}(0, k(x, x')$ (right) after sampling the domain of the \textit{Rastrigin} function ($ID = 15$) $N=50$ times.} 
    \label{fig:gp-landscape_rastrigin}
\end{subfigure}
\begin{subfigure}[b]{0.8\linewidth}
    \includegraphics[width=1\textwidth]{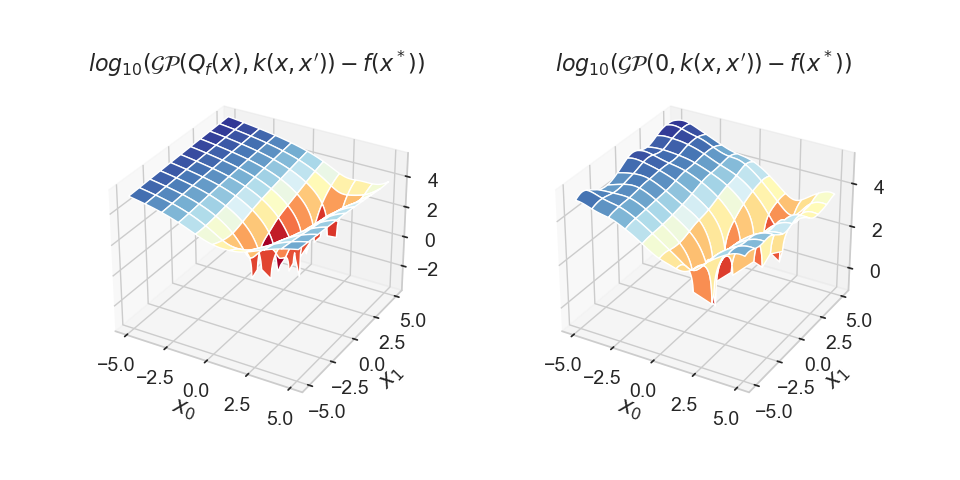}
    \caption{The prediction of the Gaussian mean after seeing the polynomial as prior $\mathop{GP}(Q_f(x), k(x, x')$ (left) compared to the prediction of a Gaussian with zero mean $\mathop{GP}(0, k(x, x')$ (right) after sampling the domain of the \textit{Schwefel} function ($ID = 20$) $N=50$ times.} 
    \label{fig:gp-landscape_schwefel}
\end{subfigure}
\caption{}
\label{fig:landscape-prediction}
\end{figure}

Finally, we investigate in Figure~\ref{fig:landscape-prediction} how the prediction of the landscape differs for the surrogate model constructed by PMBO $\mathop{GP}(Q_f(x), k(x, x'))$ and by $BO_{\mathrm{fixed}}$ $\mathop{GP}(0, k(x, x'))$.  
Figure~\ref{fig:Rastrigin-landscape} and \ref{fig:Schwefel-landscape} illustrate the ground truth of the \textit{Rastrigin} ($ID=15$) and \textit{Schwefel} ($ID=20$) function. 
The prediction of the according surrogates are shown in Figure~\ref{fig:gp-landscape_rastrigin} for the \textit{Rastrigin} and in Figure~\ref{fig:gp-landscape_schwefel} for the \textit{Schwefel} function, respectively.
The hyper-parameters are chosen as the default of the {\sc BayesO} package for the Matern$_{3/2}$ kernel.
\begin{table}[h]
\centering
\caption{This table compares the root mean square error (RMSE) of the prediction of the surrogates constructed by PMBO and $BO_{fixed}$.}
\begin{tabular}{||c c c c||} 
 \hline
 ID & name & RMSE$_{PMBO}$ & RMSE$_{BO_{\mathrm{fixed}}}$\\ [0.5ex] 
 \hline\hline
  $1$ & Sphere &  $1.65 \cdot 10^{-9}$& $1.45 \cdot 10^1$ \\ 
 \hline
  $6$ & Attractive Sector & $7.34 \cdot 10^3$ & $4.62 \cdot 10^4$\\
 \hline
  $10$ & Ellipsoidal & $2.54 \cdot 10^5$ & $6.52 \cdot 10^7$\\
 \hline
  $15$ & Rastrigin  & $1.47 \cdot 10^1$ & $3.14 \cdot 10^1$ \\
 \hline
 $20$ & Schwefel & $3.08 \cdot 10^2$ & $1.72 \cdot 10^4$ \\ 
 \hline
\end{tabular}
\label{tab:rmse}
\end{table}

Visually, the surrogate model constructed by PMBO captures the macroscopic landscape much better than $BO_\mathrm{fixed}$. 
This impression is confirmed by the prediction error for the different functions in Table~\ref{tab:rmse}.
\newpage

\newpage

\subsection{PMBO vs. other algorithms}
In the previous section, we show that the performance of PMBO is more robust towards changes of the hyper-parameters of the Gaussian process. However, the question remains whether or not this benefit makes the proposed method compatible with state-of-the-art optimization algorithms. Of special interest is the question, if PMBO can outperform Bayesian optimization that uses hyper-parameter tuning methods such as \textit{Bayesian model selection}.

We compare PMBO against the standard Bayesian optimization
(via the package {\sc BayesO} with the default settings, where the kernel is learned from the data) \cite{Kim_BayesO_A_Bayesian_2023},
CMA-ES \cite{hansen2023cma} (via the package {\sc cma} \cite{cmabbobbenchmarks}), 
and BADs \cite{acerbi2017practical} (via the package {\sc PyBads} \cite{singh2023pybads}).
\begin{figure}[h]
\centering
\begin{subfigure}[b]{0.48\linewidth}
    \includegraphics[width=1.\textwidth]
    {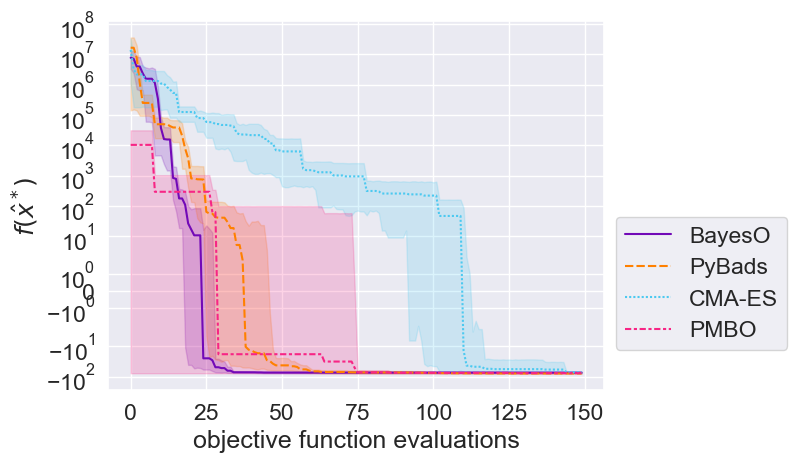}
    \caption{Ellipsoidal function}
    \label{fig:function10_d2}
\end{subfigure}
\begin{subfigure}[b]{0.48\linewidth}
    \includegraphics[width=1\textwidth]{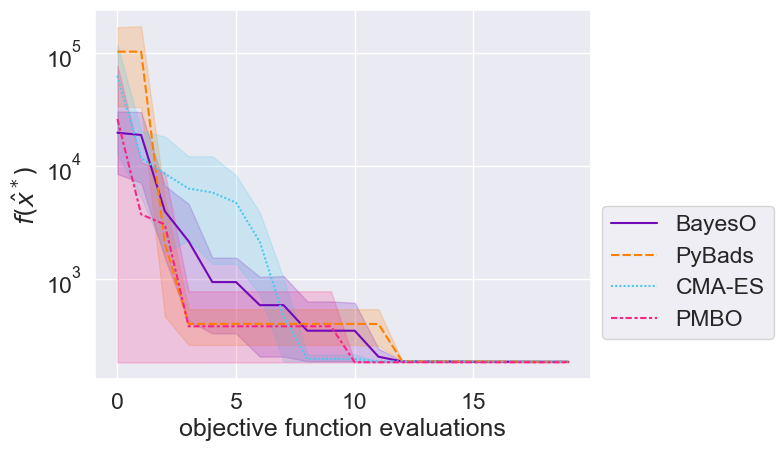}
    \caption{Schwefel function}
    \label{fig:function20_d2}
\end{subfigure}
\caption{Lowest objective function value $f(\hat{x}^*)$ per objective function evaluation for different optimization algorithms on \textit{Ellipsoidal function} $(ID=10)$ on \textit{Schwefel function} $(ID=20)$ for dimension $m=2$.}
\label{fig:optimizer-benchmark-d2}
\end{figure}

A selection of the results is shown as a representive cut through the functions in Figure~\ref{fig:optimizer-benchmark-d2} for the \textit{Ellipsoidal} function (ID=$10$) and the \textit{Schwefel} function (ID=$20$) in dimension $m=2$. 
The performance of PMBO is indeed comparable with the other chosen state-of-art algorithms.
Even though the standard deviation tends to be higher during parts of the iterations,
the average converges as fast and in some cases faster than the other algorithms. 
Especially, PMBO reaches similar or better performance in direct comparison with {\sc BayesO}, except for the function \textit{Attractive Sector} in dimension $m=5$.

\begin{figure}[h]
\centering
\begin{subfigure}[b]{0.48\linewidth}
    \includegraphics[width=1.\textwidth]
    {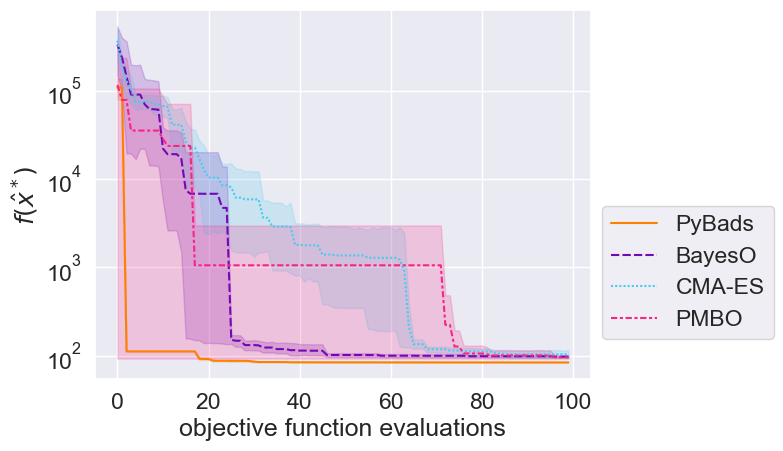}
    \caption{Attractive Sector}
    \label{fig:function6_d5}
\end{subfigure}
\begin{subfigure}[b]{0.48\linewidth}
    \includegraphics[width=1\textwidth]{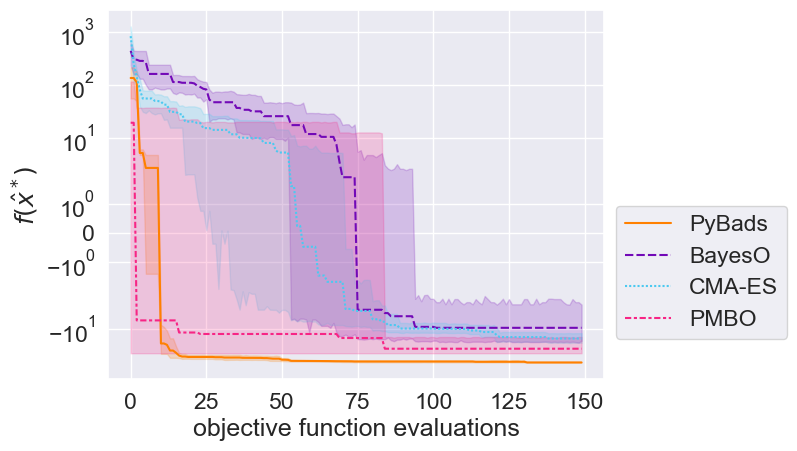}
    \caption{Rastrigin function}
    \label{fig:function15_d5}
\end{subfigure}
\caption{Lowest objective function value $f(\hat{x}^*)$ per objective function evaluation for different optimization algorithms on \textit{Attractive Sector} $(ID=6)$ on \textit{Rastrigin function} $(ID=15)$ for dimension $m=5$.}
\label{fig:optimizer-benchmark-d5}
\end{figure}

The results for dimension $m=5$ are shown in Figure~\ref{fig:optimizer-benchmark-d5}.
Here, PMBO performs poorly but still comparable with CMA-ES for the \textit{Attractive Sector} function ($ID=6$).
The opposite is the case for the \textit{Rastrigin} function. Here, PMBO outperforms the other methods. The examples showcased here are meeting the expectations of the \textit{No Free Lunch Theorem} \cite{wolpert1997no} which, loosely speaking, states that no universally best-performing algorithm exist. 

\begin{figure}[h!]
\centering
\begin{subfigure}[b]{0.495\linewidth}
    \includegraphics[width=1.\textwidth]
    {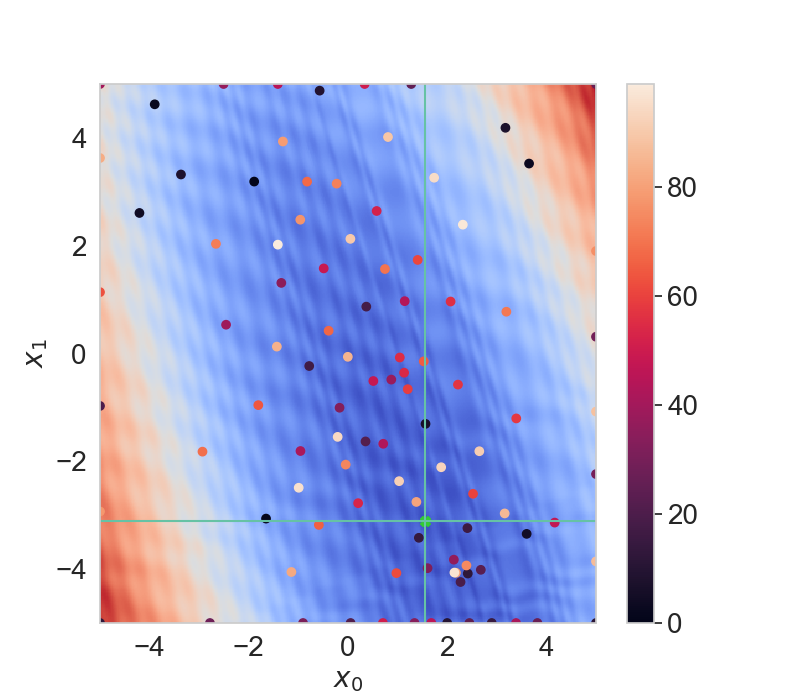}
    \caption{BayesO}
\end{subfigure}
\begin{subfigure}[b]{0.495\linewidth}
    \includegraphics[width=1\textwidth]{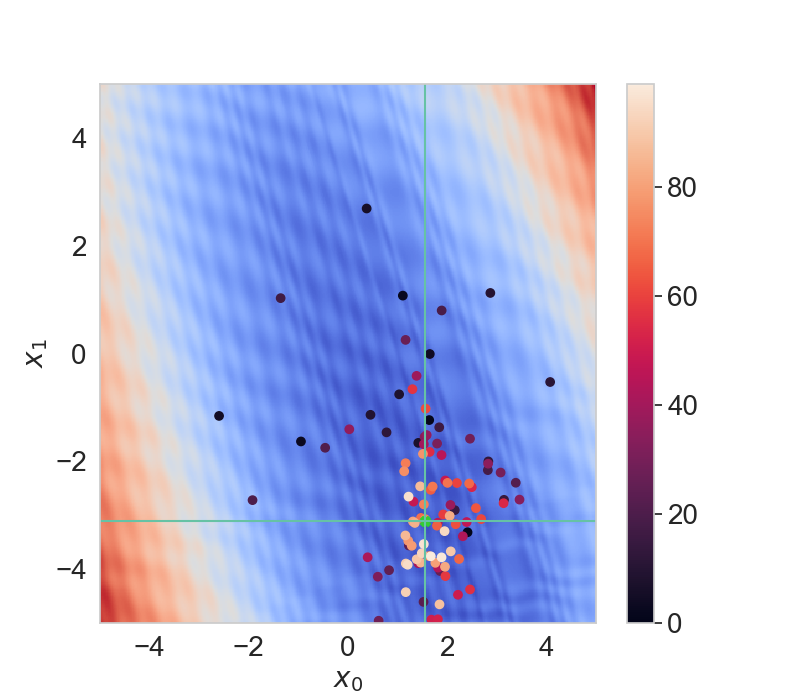}
    \caption{CMA-ES}
\end{subfigure}
\begin{subfigure}[b]{0.495\linewidth}
    \includegraphics[width=1\textwidth]{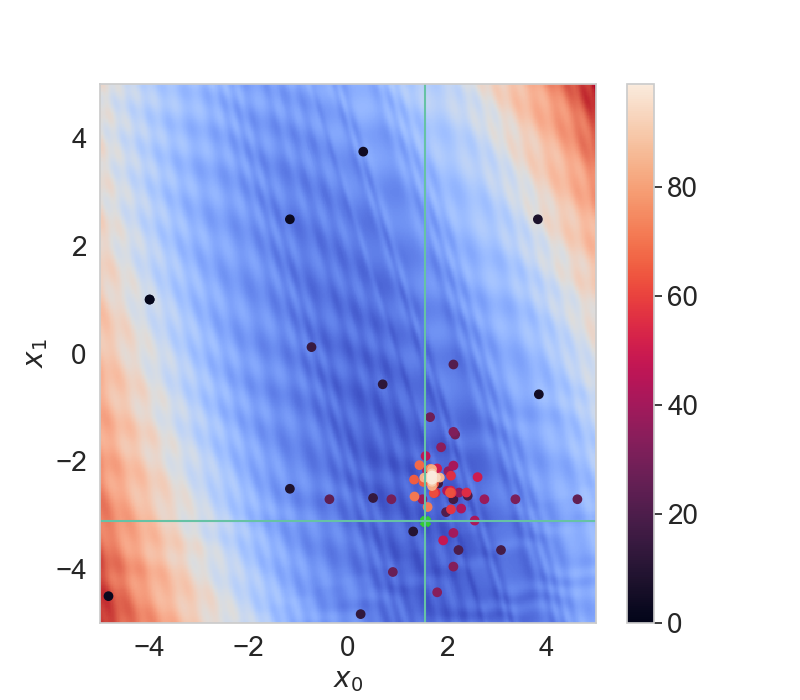}
    \caption{PyBads}
\end{subfigure}
\begin{subfigure}[b]{0.495\linewidth}
    \includegraphics[width=1\textwidth]{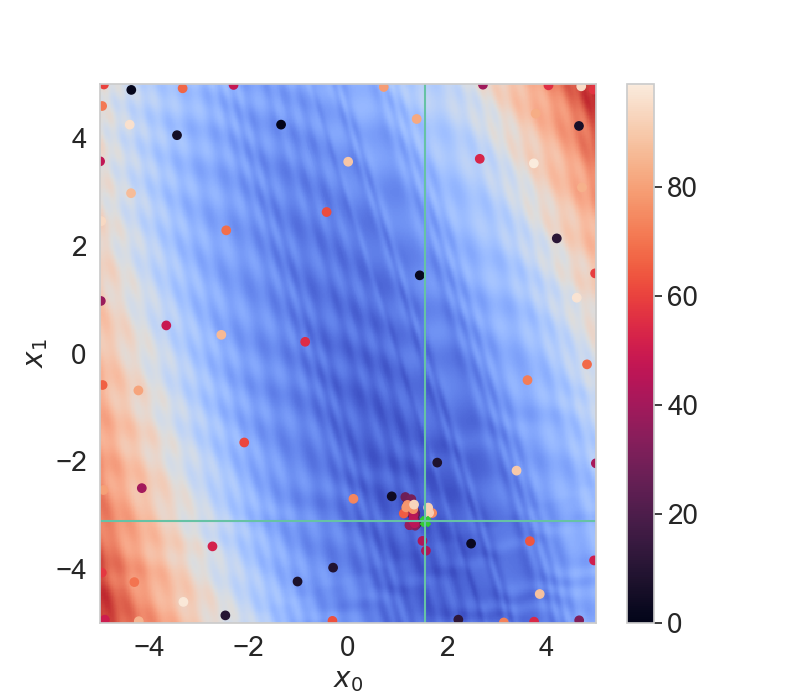}
    \caption{PMBO}
\end{subfigure}
\caption{The search behavior of {\sc BayesO}, {\sc CMA-ES}, {\sc PyBads} and {\sc PMBO} on the first $100$ samples on the \textit{Rastrigin} function ($ID = 15$) in $m=2$. The center of the green cross denotes the true minimum and the colors of the scatter plot encode the iteration number.}
\label{fig:search-behavior}
\end{figure}
Furthermore, we investigate how the different algorithms sample the search space.
Figure~\ref{fig:search-behavior} illustrates the sampling strategies for the \emph{Rastrigin} function ($ID = 15$) in dimension $m=2$. {\sc BayesO} chooses samples in the whole domain, while more samples cluster around the global minimum.
{\sc cma} starts with an initial mean and a standard deviation of a multivariate Gaussian distribution, of which the mean in the example shown here is already close to the minimum, but the algorithm eventually still fails to reach it.
It seems {\sc cma} gets stuck in a local optimum.

During the first iterations, {\sc PyBads} samples over the whole domain and then quickly concentrates samples close to the optimum.
Meanwhile, our approach takes samples across the whole domain, and like {\sc BayesO}, quite a few of these samples cover the boundaries.
With increasing sample size, {\sc PMBO} concentrates samples around the global optimum.

For the sake of completeness, it is important to mention that for the same setup but with more or less luckily chosen initial samples, all algorithms vary a bit in their performance, either for the better or for the worse.
In the given example, cases have been observed where all considered algorithms get stuck in a local optimum.

\subsection{Discussion} 
In the results above, we demonstrate the performance of PMBO relative to Bayesian optimization with fixed hyper-parameters and we compare it with several state-of-the-art algorithms.

Bayesian optimization is often the pivotal choice among model-based optimization algorithm. However, its blind spot comes by its strong sensitivity towards the choice of its hyper-parameters, whose selection is a highly non-trivial task. As shown in Section~\ref{sec:PMBOvsBO}, PMBO mitigates this sensitivity by fitting a polynomial which imposes a global structure on the objective landscape.
Together with the notion of uncertainty inherited from Bayesian optimization, PMBO opens doors for further analysis, e.g., of the extrema and near-optimal feasible regions.
These benefits are completely missing from evolutionary algorithms.

PMBO is expected to show strong performance for functions that from a macroscopic perspective resemble \textit{Bos-Levenberg-Trefethen (BLT)} functions (Section~\ref{sec:MIT}).
However, for functions with high local variations, we expect evolutionary algorithms (like CMA-ES) or hybrid solutions (like BADs) to perform better.

Until now, PMBO is limited to low-dimensional problems and works well up to dimension $m=5$.
The cause of this limitation is the poor scalability of updating the complexity of the polynomial:
despite a sub-exponential dependence of the number of interpolation nodes (coefficients) on the polynomial degree and the dimension,
the number of samples needed to accurately and stably solve the least-squares problem of Eq.~(\ref{eq:REG}) still exceeds the budget typically available in practical high-dimensional optimization problems.

Instead of accurately approximating the global objective landscape over the whole domain, approximating the landscape locally (in the region where the optima are expected to be) and with a relatively low degree polynomial (requiring a relatively small sample size to fit) may be sufficient.
However, a solution for adaptively refining the polynomial model in the regions of interest is -- to our current knowledge -- yet to be found. 

\section{Conclusion}
We proposed the Polynomial-Model-Based Optimization (PMBO) algorithm, which addresses the challenge of model-based black-box optimization by combining the best of both worlds from Gaussian processes and polynomial regression. 

We show that PMBOs' performance is largely independent from the choice of Gaussian process parameters and can successfully compete with state-of-the-art evolutionary, model-based, and hybrid algorithms. 

The PMBO-surrogate-model provides a baseline for hybridisation with evolutionary algorithms (similar to the approach of BADS),
which may further improve their performance.
Such hybridization, applications to stochastic problems and parallelization are part of future work.
\section*{Acknowledgements}
This work was partially funded by the Center of Advanced Systems Understanding (CASUS), financed by Germany's Federal Ministry of Education and Research (BMBF) and by the Saxon Ministry for Science, Culture and Tourism (SMWK) with tax funds on the basis of the budget approved by the Saxon State Parliament. \\

Pau Batlle gratefully acknowledges support from the Air Force Office of Scientific Research under MURI award number FA9550-20-1-0358 (Machine Learning and Physics-Based Modeling and Simulation), the Jet Propulsion Laboratory, California Institute of Technology, under a contract with the National Aeronautics and Space Administration, and Beyond Limits (Learning Optimal Models) through CAST (The Caltech Center for Autonomous Systems and Technologies).
\bibliographystyle{ACM-Reference-Format}
\bibliography{ref.bib}

\end{document}